\documentclass{amsart}
\usepackage{amsmath,amssymb}
\usepackage{yhmath}

\newcommand{\bdis}{\begin{displaymath}}
\newcommand{\edis}{\end{displaymath}}
\newcommand{\be}{\begin{equation}}
\newcommand{\ee}{\end{equation}}
\newcommand{\mbb}{\mathbb}
\newcommand{\mcal}{\mathcal}

\newcommand{\vp}{\varphi} 
\newcommand{\vth}{\vartheta}

\newcommand{\zf}{\zeta\left(\frac{1}{2}+it\right)}

\newcommand{\FRm}{\frac{x^m+y^m}{z^m}}

\DeclareMathOperator{\im}{Im}

\theoremstyle{definition}

\theoremstyle{remark}
\newtheorem{remark}[]{Remark}

\newtheorem*{mydef11}{{\bf Theorem 1}}

\newtheorem*{mydef12}{{\bf Theorem 2}}

\newtheorem*{mydef13}{{\bf Theorem 3}}

\newtheorem*{mydef14}{{\bf Theorem 4}}

\newtheorem*{mydef15}{{\bf Theorem 5}} 

\newtheorem*{mydef16}{{\bf Theorem 6}} 

\newtheorem*{mydef17}{{\bf Theorem 7}} 

\newtheorem*{mydef18}{{\bf Theorem 8}} 

\newtheorem*{mydef19}{{\bf Theorem 9}} 

\newtheorem*{mydef110}{{\bf Theorem 10}}

\newtheorem*{mydef41}{{\bf Corollary 1}}

\newtheorem*{mydef42}{{\bf Corollary 2}}

\newtheorem*{mydef51}{{\bf Lemma 1}}

\newtheorem*{mydef52}{{\bf Lemma 2}}

\newtheorem*{mydef81}{{\bf Property 1}}

\newtheorem*{mydef82}{{\bf Property 2}}

\numberwithin{equation}{section}

\begin{document}

\title[Jacob's ladders, new classes of sum variants \dots]{Jacob's ladders, new classes of sum variants of almost linear formula and corresponding $\zeta$-functionals and $\zeta$-equivalents of the Fermat-Wiles theorem}

\author{Jan Moser}

\address{Department of Mathematical Analysis and Numerical Mathematics, Comenius University, Mlynska Dolina M105, 842 48 Bratislava, SLOVAKIA}

\email{jan.mozer@fmph.uniba.sk}


\begin{abstract}
In this paper we obtain new classes of the second generation of new $\zeta$-functionals and corresponding $\zeta$-equivalents of the Fermat-Wiles theorem. These are based on sets of elementary rectangular $\zeta$-signals generated by the Gramm's sequence. 
\end{abstract}
\maketitle

\section{Introduction} 

\subsection{} 

In our recent paper \cite{7} we have defined the set $\{\mcal{A}_1(n)\}_{n=1}^{N(T)}$ of the rectangular signals $\mcal{A}_1(n)$ and we have obtained the following asymptotic formula 
\be \label{1.1} 
\sum_{n=1}^N |\mcal{A}_1(n)|\sim (1-c)T,\ T\to\infty. 
\ee 
This formula represents the third sum variant of our almost linear formula\footnote{See \cite{6}, (3.4), (3.6).} 
\be \label{1.2} 
\int_T^{\overset{1}{T}(T)}Z^2(t){\rm d}t\sim (1-c)T,\ \overset{1}{T}(T)=\vp_1^{-1}(T),\ T\to\infty 
\ee 
for the Riemann real-valued function 
\be \label{1.3} 
Z(t)=e^{i\vth(t)}\zf, 
\ee  
where\footnote{See \cite{8}, (35), (44), (62), comp. \cite{10}, p. 98.} 
\be \label{1.4} 
\begin{split}
	& \vth(t)=-\frac t2\ln\pi+\im\ln\Gamma\left(\frac 14+i\frac t2\right)= \\ 
	& \frac t2\ln\frac{t}{2\pi}-\frac t2-\frac{\pi}{8}+\mcal{O}\left(\frac 1t\right), 
\end{split}
\ee 
and, of course, 
\be \label{1.5} 
Z^2(t)=\left|\zf\right|^2. 
\ee 

\begin{remark}
	Let us remind that the set $\{\mcal{A}_1(n)\}$ was constructed on the base of the sequence 
	\be \label{1.6} 
	\{\gamma_k\}_{k=1}^\infty,\ \gamma_k\in (0,+\infty) 
	\ee 
	of the zeros of the function $Z(t)$. 
\end{remark} 

\subsection{} 

In present paper we use the Gramm's sequence 
\be \label{1.7} 
\{t_\nu\}_{\nu=1}^\infty:\ \vth(t_\nu)=\nu\pi 
\ee 
as the base of the construction, that is we use the sequence of roots of the equation $\vth(t)=\nu \pi$ instead of the sequence (\ref{1.6}). 

On this base we obtain the new set 
\bdis 
\{\mcal{A}_2(n)\},\ n=1,\dots, N(T)
\edis 
of rectangular signals and corresponding asymptotic formula: 
\be \label{1.8} 
\sum_{n=1}^N |\mcal{A}_2(n)|\sim (1-c)T,\ T\to\infty, 
\ee  
which represents the fourth sum variant of our almost linear formula (\ref{1.2}), and also 
\be \label{1.9} 
\frac{1}{N}\sum_{n=1}^N|\mcal{A}_2(n)|\sim 2\pi,\ T\to\infty
\ee 
that controls the oscillations of the areas $|\mcal{A}_2(n)|$ around the area (for example) 
\be \label{1.10} 
2\pi=\pi(\sqrt{2})^2
\ee 
of the circle with the radius $R=\sqrt{2}$. 

\begin{remark}
A formula similar to (\ref{1.9}) holds true conditionally also for the rectangles $\mcal{A}_1(n)$, but only on the very strong Mertens hypothesis. Consequently, there is a big difference between the results mentioned here (comp. the sixth part) and the result (\ref{1.9}) that holds true independently on any unproved hypothesis. 
\end{remark} 

\subsection{} 

Next, if we write the formula (\ref{1.8}) in the form 
\be \label{1.11} 
\sum_{T<\tilde{t}(n)<\overset{1}{T}(T)}|\mcal{A}_2(n)|\sim (1-c)T,\ T\to\infty, 
\ee 
where the symbol $\tilde{t}(n)$ denotes the mean value of the function 
\be \label{1.12} 
Z^2(t),\ t\in[t_n,t_{n+1}]\subset (T,\overset{1}{T}(T)), 
\ee  
and $t_n$ is the reduced form of the Gramm's point $t_{\nu+n}$, then we obtain, for example, the following functional 
\be \label{1.13} 
\lim_{\tau\to\infty}\frac{1}{\tau}\left\{
\sum_{\tilde{t}(n)>\frac{x}{1-c}\tau}^{\tilde{t}(n)<[\frac{x}{1-c}\tau]^1}|\mcal{A}_2(n)|
\right\} = x 
\ee   
for every fixed $x>0$, and next, in the special case of Fermat's rationals 
\be \label{1.14} 
x\to \FRm,\ x,y,z,m\in\mbb{N},\ m\geq 3, 
\ee  
the following new $\zeta$-condition 
\be \label{1.15} 
\lim_{\tau\to\infty}\frac{1}{\tau}\left\{
\sum_{\tilde{t}(n)>\FRm\frac{\tau}{1-c}}^{\tilde{t}(n)<[\FRm\frac{\tau}{1-c}]^1}|\mcal{A}_2(n)|
\right\}\not=1
\ee 
on the set of all Fermat's rationals, that expresses the new $\zeta$-equivalent of the Fermat-Wiles theorem, based on the Gramm's sequence $\{t_\nu\}$. 

\subsection{} 

Finally, we obtain for corresponding parts of the rectangles $|\mcal{A}_2(n)|$ also more complicated formulas. For the sake of simplicity we demonstrate some of these as the symbols\footnote{Complete explanations can be found in the fifth part of this work.}: 
\begin{itemize}
	\item[(a)] the set of functionals 
	\be \label{1.16} 
	\lim_{\tau\to\infty}\frac{1}{\tau}\left\{
	\sum_{t_n>\frac{x}{1-c}\tau}^{t_n<[\frac{x}{1-c}\tau]^1}\sum_{i=1}^l|\mcal{A}_2^{r_i,p_i}(n)|
	\right\}=\frac{l}{sq}x, 
	\ee 
	for every fixed 
	\bdis 
	\begin{split}
	& x>0,\ s,q\geq 2,\ 2\leq l\leq sq-1, \\ 
	& 1\leq r_i\leq s,\ 1\leq p_i\leq q, 
	\end{split}
	\edis 
	\item[(b)] the set of corresponding $\zeta$-equivalents of the Fermat-Wiles theorem 
	\be \label{1.17} 
	\lim_{\tau\to\infty}\frac{1}{\tau}\left\{
	\sum_{t_n>\FRm\frac{\tau}{1-c}}^{t_n<[\FRm\frac{\tau}{1-c}]^1}\sum_{i=1}^l|\mcal{A}_2^{r_i,p_i}(n)|
	\right\}\not= \frac{l}{sq}. 
	\ee 
\end{itemize} 

\begin{remark}
The symbols (\ref{1.16}) and (\ref{1.17}) represent the formulae of the $l$-class of the second generation for the formulae (\ref{1.13}) and (\ref{1.15}). Let us notice explicitly that all our results on functionals and corresponding $\zeta$-equivalents of the Fermat-Wiles theorem based on our almost linear formula (\ref{1.2}) represent simultaneously new results about the structure of the classical Hardy-Littlewood integral (1918)\footnote{See \cite{1}.} 
\bdis 
\int_0^T Z^2(t){\rm d}t
\edis 
since our integral, see (\ref{1.2}), 
\bdis 
\int_T^{\overset{1}{T}(T)}Z^2(t){\rm d}t 
\edis 
stands for new class of increments of the Hardy-Littlewood integral. 
\end{remark}

\section{Jacob's ladders: notions and basic geometrical properties}  

\subsection{}

In this paper we use the following notions of our works \cite{2} -- \cite{6}: 
\begin{itemize}
\item[{\tt (a)}] Jacob's ladder $\vp_1(T)$, 
\item[{\tt (b)}] direct iterations of Jacob's ladders 
\bdis 
\begin{split}
	& \vp_1^0(t)=t,\ \vp_1^1(t)=\vp_1(t),\ \vp_1^2(t)=\vp_1(\vp_1(t)),\dots , \\ 
	& \vp_1^k(t)=\vp_1(\vp_1^{k-1}(t))
\end{split}
\edis 
for every fixed natural number $k$, 
\item[{\tt (c)}] reverse iterations of Jacob's ladders 
\be \label{2.1}  
\begin{split}
	& \vp_1^{-1}(T)=\overset{1}{T},\ \vp_1^{-2}(T)=\vp_1^{-1}(\overset{1}{T})=\overset{2}{T},\dots, \\ 
	& \vp_1^{-r}(T)=\vp_1^{-1}(\overset{r-1}{T})=\overset{r}{T},\ r=1,\dots,k, 
\end{split} 
\ee   
where, for example, 
\be \label{2.2} 
\vp_1(\overset{r}{T})=\overset{r-1}{T}
\ee  
for every fixed $k\in\mbb{N}$ and every sufficiently big $T>0$. We also use the properties of the reverse iterations listed below.  
\be \label{2.3}
\overset{r}{T}-\overset{r-1}{T}\sim(1-c)\pi(\overset{r}{T});\ \pi(\overset{r}{T})\sim\frac{\overset{r}{T}}{\ln \overset{r}{T}},\ r=1,\dots,k,\ T\to\infty,  
\ee 
\be \label{2.4} 
\overset{0}{T}=T<\overset{1}{T}(T)<\overset{2}{T}(T)<\dots<\overset{k}{T}(T), 
\ee 
and 
\be \label{2.5} 
T\sim \overset{1}{T}\sim \overset{2}{T}\sim \dots\sim \overset{k}{T},\ T\to\infty.   
\ee  
\end{itemize} 

\begin{remark}
	The asymptotic behaviour of the points 
	\bdis 
	\{T,\overset{1}{T},\dots,\overset{k}{T}\}
	\edis  
	is as follows: at $T\to\infty$ these points recede unboundedly each from other and all together are receding to infinity. Hence, the set of these points behaves at $T\to\infty$ as one-dimensional Friedmann-Hubble expanding Universe. 
\end{remark}  

\subsection{} 

Let us remind that we have proved\footnote{See \cite{5}.} the existence of almost linear increments 
\be \label{2.6} 
\begin{split}
& \int_{\overset{r-1}{T}}^{\overset{r}{T}}\left|\zf\right|^2{\rm d}t\sim (1-c)\overset{r-1}{T}, \\ 
& r=1,\dots,k,\ T\to\infty,\ \overset{r}{T}=\overset{r}{T}(T)=\vp_1^{-r}(T)
\end{split} 
\ee 
for the Hardy-Littlewood integral (1918), \cite{1}: 
\be \label{2.7} 
J(T)=\int_0^T\left|\zf\right|^2{\rm d}t. 
\ee  

For completeness, we give here some basic geometrical properties related to Jacob's ladders. These are generated by the sequence 
\be \label{2.8} 
T\to \left\{\overset{r}{T}(T)\right\}_{r=1}^k
\ee 
of reverse iterations of the Jacob's ladders for every sufficiently big $T>0$ and every fixed $k\in\mbb{N}$. 

\begin{mydef81}
The sequence (\ref{2.8}) defines a partition of the segment $[T,\overset{k}{T}]$ as follows 
\be \label{2.9} 
|[T,\overset{k}{T}]|=\sum_{r=1}^k|[\overset{r-1}{T},\overset{r}{T}]|
\ee 
on the asymptotically equidistant parts 
\be \label{2.10} 
\begin{split}
& \overset{r}{T}-\overset{r-1}{T}\sim \overset{r+1}{T}-\overset{r}{T}, \\ 
& r=1,\dots,k-1,\ T\to\infty. 
\end{split}
\ee 
\end{mydef81} 

\begin{mydef82}
Simultaneously with the Property 1, the sequence (\ref{2.8}) defines the partition of the integral 
\be \label{2.11} 
\int_T^{\overset{k}{T}}\left|\zf\right|^2{\rm d}t
\ee 
into the parts 
\be \label{2.12} 
\int_T^{\overset{k}{T}}\left|\zf\right|^2{\rm d}t=\sum_{r=1}^k\int_{\overset{r-1}{T}}^{\overset{r}{T}}\left|\zf\right|^2{\rm d}t, 
\ee 
that are asymptotically equal 
\be \label{2.13} 
\int_{\overset{r-1}{T}}^{\overset{r}{T}}\left|\zf\right|^2{\rm d}t\sim \int_{\overset{r}{T}}^{\overset{r+1}{T}}\left|\zf\right|^2{\rm d}t,\ T\to\infty. 
\ee 
\end{mydef82} 

It is clear, that (\ref{2.10}) follows from (\ref{2.3}) and (\ref{2.5}) since 
\be \label{2.14} 
\overset{r}{T}-\overset{r-1}{T}\sim (1-c)\frac{\overset{r}{T}}{\ln \overset{r}{T}}\sim (1-c)\frac{T}{\ln T},\ r=1,\dots,k, 
\ee  
while our eq. (\ref{2.13}) follows from (\ref{2.6}) and (\ref{2.5}).  

\section{The next sum variant of almost linear formula based on Gramm's sequence} 

\subsection{}  

We have obtained the following in our paper \cite{7}, (3.14), (3.17): let the symbol 
\be \label{3.1} 
\{\gamma_n\}_{n=0}^{N+2}
\ee 
denote the subset of the set of zeros of the function 
\be \label{3.2} 
Z(t),\ t\in (T,\overset{1}{T}(T)),\ \overset{1}{T}(T)=\vp_1^{-1}(T), 
\ee  
where 
\be \label{3.3} 
\gamma_0\leq T<\gamma_1<\gamma_2<\dots<\gamma_N<\gamma_{N+1}<\overset{1}{T}(T)\leq \gamma_{N+2},\ N=N(T), 
\ee  
i.e. the subset (\ref{3.1}) of zeros is selected from the set (\ref{1.6}) by the segment $[T,\overset{1}{T}(T)]$. Then it is true\footnote{See \cite{7}, (3.16), (3.19).}: 
\be \label{3.4} 
\sum_{T\leq \gamma_n<\overset{1}{T}(T)}|\mcal{A}_1(n)|\sim (1-c)T,\ T\to\infty , 
\ee  
or, in the form 
\be \label{3.5} 
\sum_{n=1}^{N(T)}|\mcal{A}_1(n)|\sim (1-c)T,\ T\to\infty, 
\ee  
for the set of rectangular $\zeta$-signals $\mcal{A}_1(n)$, where 
\be \label{3.6} 
\begin{split}
& \mcal{A}_1(n)=[\gamma_n,\gamma_{n+1}]\times [0,Z^2(\bar{t}(n))],\ n=1,\dots,N(T), \\ 
& Z^2(\bar{t}(n))=\frac{1}{\gamma_{n+1}-\gamma_n}\int_{\gamma_n}^{\gamma_{n+1}}Z^2(t){\rm d}t, \\ 
& |\mcal{A}_1(n)|=Z^2(\bar{t}(n))(\gamma_{n+1}-\gamma_n). 
\end{split}
\ee  

\subsection{} 

In this paper we use the Gramm' sequence as a base. Let us remind that for this sequence the following Titchmarsh formula holds true 
\be \label{3.7} 
t_{\nu+1}-t_\nu=\frac{2\pi}{\ln t_\nu}+\frac{2\pi\ln2\pi}{\ln^2t_\nu}+\mcal{O}\left(\frac{1}{\ln^3t_\nu}\right),\ t_\nu\to\infty. 
\ee 
Now, we select from the sequence (\ref{1.7}) a local sequence 
\be \label{3.8} 
\{t_{n}\}_{n=0}^{N+2},\ t_{\nu+n}\to t_n,\ N=N(T), 
\ee  
i.e. the letter $n$ denotes the local index, and 
\be \label{3.9} 
t_0\leq T<t_1<t_2<\dots<t_N<t_{N+1}<\overset{1}{T}(T)\leq t_{N+2} 
\ee 
is the reduction of (\ref{3.8}) to the interior of the segment $[T,\overset{1}{T}(T)]$.  

\subsection{} 

We have the following formula, see (\ref{3.9}): 
\be \label{3.10} 
\begin{split}
& \int_T^{\overset{1}{T}(T)}Z^2(t){\rm d}t=\sum_{n=1}^{N(T)}\int_{t_n}^{t_{n+1}}Z^2(t){\rm d}t+ \\ 
& \int_T^{t_1}Z^2(t){\rm d}t+\int_{t_{N+1}}^{\overset{1}{T}(T)}Z^2(t){\rm d}t. 
\end{split} 
\ee  
Now, it is sufficient to use the estimates\footnote{Comp. (\ref{3.7})} 
\be \label{3.11} 
\begin{split}
& t_1-T,\ \overset{1}{T}(T)-t_{N+1}=\mcal{O}\left(\frac{1}{\ln T}\right),  \\ 
& Z^2(t)=\mcal{O}(T^{1/3}),\ t\in [T,t_1]\cup [t_{N+1},\overset{1}{T}(T)], 
\end{split}
\ee 
and these imply the following result 
\be \label{3.12} 
\begin{split}
& \int_T^{\overset{1}{T}(T)}Z^2(t){\rm d}t=\sum_{n=1}^N\int_{t_n}^{t_{n+1}}Z^2(t){\rm d}t+\mcal{O}\left(\frac{T^{1/3}}{\ln T}\right)= \\ 
& \sum_{T<t_n<\overset{1}{T}(T)}\int_{t_{n}}^{t_{n+1}}Z^2(t){\rm d}t+\mcal{O}(T^{1/3+\delta}),\ T\to\infty  
\end{split}
\ee 
(the remainder in the formula (\ref{1.2}) is $\mcal{O}(T^{1/3+\delta})$, $0<\delta$ is any small fixed value). 

Consequently, the following lemma holds true (see (\ref{1.2}), (\ref{3.12})) 
\begin{mydef51}
\be \label{3.13} 
\sum_{T<t_n<\overset{1}{T}(T)}\int_{t_{n}}^{t_{n+1}}Z^2(t){\rm d}t\sim (1-c)T,\ T\to\infty, 
\ee  
or 
\be \label{3.14} 
\sum_{n=1}^{N(T)}\int_{t_n}^{t_{n+1}}Z^2(t){\rm d}t\sim (1-c)T,\ T\to\infty. 
\ee 
\end{mydef51} 

\subsection{} 

Let 
\be \label{3.15} 
Z^2(\tilde{t}(n))=\frac{1}{t_{n+1}-t_n}\int_{t_n}^{t_{n+1}}Z^2(t){\rm d}t,\ \tilde{t}_n\in (t_n,t_{n+1}), 
\ee  
then we can define the following rectangular $\zeta$-signals 
\be \label{3.16} 
\mcal{A}_2(n)=[t_n,t_{n+1}]\times [0,Z^2(\tilde{t}(n))],\ n=1,\dots,N(T)
\ee  
with the area 
\be \label{3.17} 
|\mcal{A}_2(n)|=Z^2(\tilde{t}(n))(t_{n+1}-t_n). 
\ee 
Now, the formula (\ref{3.13}) (for example) implies the result: 
\begin{mydef11}
\be \label{3.18} 
\sum_{T<t_n<\overset{1}{T}(T)}|\mcal{A}_2(n)|\sim (1-c)T,\ T\to\infty. 
\ee 
\end{mydef11} 

\begin{remark}
The formula (\ref{3.18}) represents the fourth sum variant of our almost linear formula (\ref{1.2}). 
\end{remark} 

\subsection{} 

Of course, the formula (\ref{3.18}) generates also new variants of $\zeta$-functional and $\zeta$-equivalent of the Fermat-Wiles theorem. Namely, if we use in (\ref{3.18}) the substitution 
\be \label{3.19} 
T=\frac{x}{1-c}\tau,\ \{T\to+\infty\}\Leftrightarrow\{\tau\to+\infty\}
\ee 
for every fixed $x>0$, then we obtain the following $\zeta$-functional. 

\begin{mydef12}
\be \label{3.20} 
\lim_{\tau\to\infty}\frac{1}{\tau}\left\{
\sum_{\tilde{t}(n)>\frac{x}{1-c}\tau}^{\tilde{t}(n)<[\frac{x}{1-c}\tau]^1}|\mcal{A}_2(n)|
\right\}=x;\ [G]^1=\vp_1^{-1}(G), 
\ee 
for every fixed $x>0$. 
\end{mydef12} 

And, in the special case of Fermat's rationals, see (\ref{1.14}), we obtain the following. 

\begin{mydef41}
\be \label{3.21} 
\lim_{\tau\to\infty}\frac{1}{\tau}\left\{
\sum_{\tilde{t}(n)>\FRm\frac{\tau}{1-c}}^{\tilde{t}(n)<[\FRm\frac{\tau}{1-c}]^1}|\mcal{A}_2(n)|
\right\}=\FRm
\ee 
for every fixed Fermat's rational. 
\end{mydef41} 

Consequently, the following theorem holds true. 

\begin{mydef13}
The $\zeta$-condition 
\be \label{3.22} 
\lim_{\tau\to\infty}\frac{1}{\tau}\left\{
\sum_{\tilde{t}(n)>\FRm\frac{\tau}{1-c}}^{\tilde{t}(n)<[\FRm\frac{\tau}{1-c}]^1}|\mcal{A}_2(n)|
\right\}\not=1 
\ee  
on the set of all Fermat's rationals expresses the new type of $\zeta$-equivalent of the Fermat-Wiles theorem based on the Gramm' sequence. 
\end{mydef13} 

\section{The second generation of almost linear formula} 

\subsection{} 

Here we give some partitions of the rectangles $\mcal{A}_2(n)|$, see (\ref{3.16}), on corresponding sets of elementary parts - also rectangles - as follows. 

First, we give the equidistant partition of the segments 
\be \label{4.1} 
[t_n,t_{n+1}],\ n=1,\dots,N(T)
\ee 
by means of the set of points 
\be \label{4.2} 
\{t_n^1,t_n^2,\dots,t_n^{s-1}\}
\ee 
for every fixed $n$ as follows 
\be \label{4.3} 
\begin{split}
& t_n^0=t_n<t_n^1<t_n^2<\dots<t_n^{r-1}<t_n^r<\dots<t_n^{s-1}<t_n^s=t_{n+1}, \\ 
& t_n^r-t_n^{r-1}=\frac{1}{s}(t_{n+1}-t_n),\ r=1,\dots,s. 
\end{split}
\ee 
And next, similarly to (\ref{4.1}) -- (\ref{4.3}), we have also the equidistant partition of the segment 
\be \label{4.4} 
[0,Z^2(\tilde{t}(n))]
\ee 
as follows 
\be \label{4.5} 
\begin{split}
& y_n^0=0<y_n^1<y_n^2<\dots<y_n^{p-1}<y_n^p<\dots<y_n^{q-1}<y_n^q=Z^2(\tilde{t}(n)), \\ 
& y_n^p-y_n^{p-1}=\frac{1}{q}Z^2(\tilde{t}(n)),\ p=1,\dots,q. 
\end{split}
\ee 
Now, we put\footnote{See (\ref{4.3}) and (\ref{4.5}).} 
\be \label{4.6} 
\begin{split}
& \mcal{A}_2^{r,p}(n)=[t_n^{r-1},t_n^r]\times[y_n^{p-1},y_n^p], \\ 
& r=1,\dots,s,\  p=1,\dots, q 
\end{split}
\ee 
for every fixed $n=1,\dots,N(T)$, where\footnote{See (\ref{3.17}).} 
\be \label{4.7} 
|\mcal{A}_2^{r,p}(n)|=\frac{1}{sq}(t_{n+1}-t_n)Z^2(\tilde{t}(n))=\frac{1}{sq}|\mcal{A}_2(n)|, 
\ee  
and, of course, 
\be \label{4.8} 
\mcal{A}_2(n)=\bigcup_{r=1}^s\bigcup_{p=1}^q\mcal{A}_2^{r,p}(n), 
\ee  
and 
\be \label{4.9} 
|\mcal{A}_2(n)|=\sum_{r=1}^s\sum_{p=1}^q|\mcal{A}_2^{r,p}(n)|,\ n=1,\dots,N(T). 
\ee 

\subsection{} 

The simple choice\footnote{Comp. (\ref{4.6}).} 
\be \label{4.10} 
\mcal{A}_2^{s,q}(n)\in \left\{\mcal{A}_2^{r,p}(n)\right\}_{r=1,p=1}^{s,q}
\ee 
is completely sufficient for our purpose. In fact, the elements of the set in (\ref{4.6}) are congruent rectangles for every fixed $n$. Consequently, 
\be \label{4.11}  
\begin{split} 
& |\mcal{A}_2^{s,q}(n)|=|\mcal{A}_2^{r,p}(n)|, \\ 
& r=1,\dots,s,\ p=1,\dots,q, 
\end{split} 
\ee 
and, of course, 
\be \label{4.12} 
\sum_{T<t_n<\overset{1}{T}(T)}|\mcal{A}_2^{r,p}(n)|=\sum_{T<t_n<\overset{1}{T}(T)}|\mcal{A}_2^{s,q}(n)|
\ee 
for every element of the set in (\ref{4.10}). 

\begin{remark}
The motivation for the choice in (\ref{4.10}) is as follows: 
\begin{itemize}
	\item[(a)] We wish to consider especially 
	\be \label{4.13} 
	|\mcal{A}_2^{s,q}(n)| 
	\ee 
	as the hover quanta (portions) of $Z^2$-energy in contrast with the quanta $|\mcal{A}_2^{r,1}(n)|$ on a \emph{ground}, for example. 
	\item[(b)] This wish is related with one of the ideas of Nicola Tesla concerning wireless transportations of electric power in atmosphere. 
\end{itemize}
\end{remark} 

\subsection{} 

Since\footnote{See (\ref{4.7}).} 
\be \label{4.14} 
|\mcal{A}_2^{s,q}(n)|=\frac{1}{sq}|\mcal{A}_2(n)|, 
\ee  
then by (\ref{3.18}) one obtains the following. 

\begin{mydef14}
\be \label{4.15} 
\sum_{T<t_n<\overset{1}{T}(T)}|\mcal{A}_2^{s,q}(n)|\sim \frac{1-c}{sq}T,\ T\to\infty, 
\ee  
for every fixed $s$ and $q$; $s,q\geq 2$. 
\end{mydef14} 

\section{The second generation of functionals and $\zeta$-equivalents of the Fermat-Wiles theorem} 

\subsection{} 

Now, if we use the substitution (\ref{3.19}), then we obtain the following functional. 

\begin{mydef15}
\be \label{5.1} 
\lim_{\tau\to\infty}\frac{1}{\tau}\left\{
\sum_{t_n>\frac{x}{1-c}\tau}^{t_n<[\frac{x}{1-c}\tau]^1}|\mcal{A}_2^{s,q}(n)|
\right\}=\frac{1}{sq}x 
\ee  
for every fixed $x>0$ and $s,q\geq 2$. 
\end{mydef15}

\begin{remark}
The formula (\ref{5.1}) represents one element of the first class of functionals of the second generation of our first functional\footnote{See \cite{6}, (4.5).} 
\be \label{5.2} 
\lim_{\tau\to\infty}\frac{1}{\tau}\int_{\frac{x}{1-c}\tau}^{[\frac{x}{1-c}\tau]^1}\left|\zf\right|^2{\rm d}t=x,\ x>0. 
\ee 
\end{remark} 

\subsection{} 

In the special case of the Fermat's rationals (\ref{1.14}) we obtain the formula: 

\begin{mydef42}
\be \label{5.3} 
\lim_{\tau\to\infty}\frac{1}{\tau}\left\{
\sum_{t_n>\FRm\frac{\tau}{1-c}}^{t_n<[\FRm\frac{\tau}{1-c}]^1}|\mcal{A}_2^{s,q}(n)|
\right\}=\frac{1}{sq}\FRm
\ee 
for every fixed Fermat's rational and every fixed $s,q\geq 2$. 
\end{mydef42} 

Consequently, the following theorem holds true. 

\begin{mydef16}
The $\zeta$-condition 
\be \label{5.4} 
\lim_{\tau\to\infty}\frac{1}{\tau}\left\{
\sum_{t_n>\FRm\frac{\tau}{1-c}}^{t_n<[\FRm\frac{\tau}{1-c}]^1}|\mcal{A}_2^{s,q}(n)|
\right\}\not=\frac{1}{sq}
\ee 
on the set of all Fermat's rationals and every fixed $s,q\geq 2$ expresses the new type of $\zeta$-equivalent of the Fermat-Wiles theorem based on the Gramm' sequence. 
\end{mydef16} 

\begin{remark}
Of course, corresponding remarks similar to the Remark 7, hold true also for the formulae (\ref{5.3}) and (\ref{5.4}). 
\end{remark} 

\subsection{} 

Now, we give some remarks about the next classes of formulas similar to (\ref{5.1}) and (\ref{5.4}).  

First, we choose the set of sequences 
\be \label{5.5} 
\begin{split}
& \{(r_i(n),p_i(n))\}_{i=1}^l,\ 2\leq l\leq sq-1, \\ 
& 1\leq r_i(n)\leq s,\ 1\leq p_i(n)\leq q 
\end{split}
\ee 
for every fixed $l$ and $n$. Next, we reduce\footnote{Comp. (\ref{4.10}), (\ref{4.11}).} this selection for every fixed $l$ on the following one 
\be \label{5.6} 
\{(r_i,p_i)\}_{i=1}^l,\ n=1,\dots,N(T). 
\ee  
Consequently, there is the functional according to the next theorem that corresponds to the choice (\ref{5.6}), and, of course, also corresponding $\zeta$-equivalent of the Fermat-Wiles theorem, see (\ref{5.1}) and (\ref{5.4}). 

\begin{mydef17}
\be \label{5.7} 
\lim_{\tau\to\infty}\frac{1}{\tau}\left\{
\sum_{t_n>\frac{x}{1-c}\tau}^{t_n<[\frac{x}{1-c}\tau]^1}\sum_{i=1}^l|\mcal{A}_2^{r_i,p_i}(n)|
\right\}=\frac{l}{sq}x 
\ee 
for every fixed 
\bdis 
x>0,\ s,q\geq 2,\ 2\leq l\leq sq-1. 
\edis 
\end{mydef17} 

\begin{remark}
If we put $l=sq$ in eq. (\ref{5.7}), then we obtain the formula\footnote{Comp. (\ref{4.8}).}: 
\bdis  
\lim_{\tau\to\infty}\frac{1}{\tau}\left\{
\sum_{t_n>\frac{x}{1-c}\tau}^{t_n<[\frac{x}{1-c}\tau]^1}|\mcal{A}_2(n)|
\right\}=x.  
\edis 
\end{remark} 

\begin{mydef18}
The $\zeta$-condition 
\be \label{5.8} 
\lim_{\tau\to\infty}\frac{1}{\tau}\left\{
\sum_{t_n>\FRm\frac{\tau}{1-c}}^{t_n<[\FRm\frac{\tau}{1-c}]^1}\sum_{i=1}^l|\mcal{A}_2^{r_i,p_i}(n)|
\right\}\not=\frac{l}{sq} 
\ee 
on the set of all Fermat's rationals and every fixed 
\bdis 
s,q\geq 2,\ 2\leq l\leq sq-1
\edis  
expresses the new type of $\zeta$-equivalent of the Fermat-Wiles theorem. 
\end{mydef18} 

\begin{remark}
Each of the formulae (\ref{5.7}) and (\ref{5.8}) represents one element of the $l$-th class of the second generation of the corresponding functionals and $\zeta$-equivalents of the Fermat-Wiles theorem. 
\end{remark} 

\section{On the canonical formula for the arithmetic mean of the set of areas $|\mcal{A}_2(n)|$} 

\subsection{} 

Since\footnote{See (\ref{2.3}), (\ref{2.5}).} 
\be \label{6.1} 
\overset{1}{T}(T)-T\sim (1-c)\frac{T}{\ln T},\ T\to\infty, 
\ee  
then 
\be \label{6.2} 
t_\nu\in [T,\overset{1}{T}(T)] \ \Rightarrow \ \frac{1}{\ln t_\nu}=\frac{1}{\ln T}+\mcal{O}\left(\frac{1}{\ln^3T}\right), 
\ee  
and, consequently, it follows from the Titchmarsh formula (\ref{3.7}) 
\be \label{6.3} 
t_{\nu+1}-t_\nu=\{1+o(1)\}\frac{2\pi}{\ln T},\ T\to\infty. 
\ee  

\subsection{} 

Next\footnote{Comp. (\ref{3.16}).}, 
\bdis 
\{\mcal{A}_2(n)\}_{n=1}^{N(T)}, 
\edis   
where 
\be \label{6.4} 
N=\sum_{(t_n,t_{n+1})\subset (T,\overset{1}{T}(T))}1, 
\ee  
and 
\be \label{6.5} 
t_{n+1}-t_n=\{1+o(1)\}\frac{2\pi}{\ln T},\ n=1,\dots,N(T),\ T\to\infty. 
\ee 
Since, see (\ref{6.1}), 
\be \label{6.6} 
\overset{1}{T}(T)-T=\{1+\bar{o}(1)\}(1-c)\frac{T}{\ln T}, 
\ee  
then it is true that 
\be \label{6.7} 
N<\frac{\overset{1}{T}(T)-T}{\{1-|o(1)|\}\frac{2\pi}{\ln T}}<\frac{1+|\bar{o}(1)|}{1-|o(1)|}\frac{1-c}{2\pi}T,\ T\to\infty, 
\ee  
and, by a similar way 
\be \label{6.8} 
N>\frac{1-|\bar{o}(1)|}{1+|o(1)|}\frac{1-c}{2\pi}T,\ T\to\infty. 
\ee 
Consequently, it follows from the last two formulae that 
\be \label{6.9} 
\lim_{\tau\to\infty}\frac{N(T)}{\frac{1-c}{2\pi}T}=1, 
\ee 
i.e. the following lemma holds true. 

\begin{mydef52}
\be \label{6.10} 
N(T)\sim \frac{1-c}{2\pi}T,\ T\to\infty. 
\ee 
\end{mydef52} 

\subsection{} 

Now, if we rewrite the formula (\ref{3.18}) in the form\footnote{Comp. (\ref{3.13}), (\ref{3.14}).}
\be \label{6.11} 
\sum_{n=1}^{N(T)}|\mcal{A}_2(n)|\sim (1-c)T,\ T\to\infty, 
\ee 
then the quotient of the formulae (\ref{6.11}) and (\ref{6.10}) implies the following result. 

\begin{mydef19}
\be \label{6.12} 
\frac{1}{N}\sum_{n=1}^{N(T)}|\mcal{A}_2(n)|\sim 2\pi,\ T\to\infty. 
\ee 
\end{mydef19} 

\begin{remark}
The asymptotic formula (\ref{6.12}) can be interpreted as: The areas $|\mcal{A}_2(n)|$ of the rectangular signals $\mcal{A}_2(n)$ generated by the Gramm' sequence $\{t_\nu\}$ are oscillating about the value $2\pi$. Next, the value $2\pi$ itself can be interpreted as: 
\begin{itemize}
	\item[(a)] In the case 
	\be \label{6.13} 
	2\pi = \pi (\sqrt{2})^2
	\ee 
	this value stands for the area of circle with radius $R=\sqrt{2}$, where $\sqrt{2}$ is the first \emph{Greek} irrational that is pleasant from the aesthetic point of view. 
	\item[(b)] In the case 
	\be \label{6.14} 
	\mcal{R}=[1,0]\times [0,2\pi];\ |\mcal{R}|=2\pi, 
	\ee  
	the value $2\pi$ expresses the area of the rectangle $\mcal{R}$ that is more corresponding to the given set $\{\mcal{A}_2(n)\}$ of the rectangular signals. 
\end{itemize}
\end{remark} 

\subsection{} 

Let us remind that there is a big divergence between the result (\ref{6.12}) based upon the Gramm' sequence $\{t_\nu\}$ and the analogical result (which we give below) for the sequence $\{\gamma_k\}$ i.e. for the set of rectangular signals $\{\mcal{A}_1(n)\}$. 

Namely, the result (\ref{6.12}) holds true independently on any unproved hypothesis, while the analogical result for the set $\{\mcal{A}_1(n)\}$ can be given only on the base of very strong Mertens hypothesis i.e. the hypothesis that 
\be \label{6.15} 
M(x)=\sum_{n\leq x}\mu(n)=\mcal{O}(\sqrt{x}),\ x\geq 1, 
\ee  
where 
\be \label{6.16} 
\frac{1}{\zeta(s)}=\sum_{n=1}^\infty \frac{\mu(n)}{n^s},\ s=\sigma+it,\ \sigma>1, 
\ee 
and $\mu(n)$ is the M\" obius function. 

It follows from the Mertens hypothesis, see \cite{7}, pp. 78, 79, concerning Cramer's and E. Landau's classical results: 
\begin{itemize}
	\item[(a)] The validity of the Riemann hypothesis. 
	\item[(b)] Simplicity of all the zeros of the Riemann zeta-function. 
\end{itemize} 

Now, it follows from (b) in our case\footnote{Comp. (\ref{6.10}).}: 
\be \label{6.17} 
N=\sum_{T<\gamma_n<\overset{1}{T}(T)}1\sim\frac{1}{2\pi}(\overset{1}{T}(T)-T)\ln T \sim \frac{1-c}{\pi}T,\ T\to\infty. 
\ee 

Consequently, we have, see (\ref{1.1}) and (\ref{6.17}), the following conditional result for the set $\mcal{A}_1(n)$ of the rectangular signals generated by the sequence $\{\gamma_k\}$: 

\begin{mydef110} 
On Mertens hypothesis the formula 
\be \label{6.18} 
\frac{1}{N}\sum_{n=1}^{N(T)}|\mcal{A}_1(n)|\sim 2\pi,\ T\to\infty. 
\ee 
holds true. 
\end{mydef110}

I would like to thank Michal Demetrian for his moral support of my study of Jacob's ladders.

\end{document}